\documentclass[12pt]{article}
\title{Induced nilpotent orbits and birational geometry}
\author{Yoshinori Namikawa}
\date{ }
\usepackage{amsmath,amssymb,amsthm, amscd}
\chardef\bslash=`\\

\def\0{{\mathcal O}}
\def\g{{\mathfrak g}}
\def\q{{\mathfrak q}}
\def\p{{\mathfrak p}}
\def\h{{\mathfrak h}}

\def\l{{\mathfrak l}}
\def\m{{\mathfrak m}}

\def\n{{\mathfrak n}}

\begin{document}
\maketitle

\begin{quote}
\hspace{2.0cm} To the memory of  Professor Masayoshi Nagata 
\end{quote}

\begin{center} 
{\bf Introduction}
\end{center}

Let $G$ be a complex simple algebraic group and let 
$\g$ be its Lie algebra. A nilpotent orbit $\mathcal{O}$ in $\g$ 
is an orbit of a nilpotent element of $\g$ by the adjoint action 
of $G$ on $\g$. Then $\mathcal{O}$ admits a natural symplectic 
2-form $\omega$ and the nilpotent orbit closure $\bar{\mathcal O}$ has symplectic  
singularities in the sense of [Be] and [Na3] (cf. [Pa], [Hi]).  
In [Ri], Richardson introduced the notion of 
so-called the {\em Richardson orbit}. 
A nilpotent orbit $\mathcal{O}$ is called Richardson if there is 
a parabolic subgroup $Q$ of $G$ such that $\mathcal{O} \cap n(\q)$ 
is an open dense subset of $n(\q)$, where $n(\q)$ is the nil-radical 
of $\q$. Later, Lusztig and Spaltenstein [L-S] generalized this notion 
to the {\em induced orbit}. A nilpotent orbit $\mathcal{O}$ is an 
induced orbit if there are a parabolic subgroup $Q$ of $G$ and 
a nilpotent orbit $\mathcal{O}'$ in the Levi subalgebra $\l(\q)$ of 
$\q := \mathrm{Lie}(Q)$ such that $\mathcal{O}$ meets $n(\q) + \mathcal{O}'$ 
in an open dense subset. If $\mathcal{O}$ is an induced orbit, one has 
a natural map (cf. (1.2)) 
$$ \nu: G \times^Q (n(\q) + \bar{\mathcal{O}}') \to \bar{\mathcal{O}}. $$
The map $\nu$ is a generically finite, projective, surjective map.   
This map is called the {\em generalized Springer map}. 
In this paper, we shall study the induced orbits from the view point 
of {\em birational geometry}. For a Richardson orbit $\mathcal{O}$, 
the Springer map $\nu$ is a map from the cotangent bundle $T^*(G/Q)$ 
of the flag variety $G/Q$ to $\bar{\mathcal O}$.  In [Fu], Fu proved 
that, if $\bar{\mathcal{O}}$ has a crepant (projective) resolution, it 
is a Springer map. Note that $Q$ is 
not unique (even up to the conjugate) for a Richardson orbit $\mathcal{O}$. 
This means that $\bar{\mathcal{O}}$ has many different crepant 
resolutions. In [Na], the author has given a description of all crepant resolutions 
of $\bar{\mathcal{O}}$ and proved that any two different 
crepant resolutions are connected by {\em Mukai flops}. 
The purpose of this paper is to generalize these to {\em all} nilpotent 
orbits $\mathcal{O}$. If $\mathcal{O}$ is not Richardson, $\bar{\mathcal{O}}$ 
has no crepant resolution. The substitute of a crepant resolution, is a 
{\em {\bf Q}-factorial terminalization}.           
Let $X$ be a complex algebraic variety with rational Gorenstein 
singularities. 
A partial resolution $f: Y \to X$ of $X$ is said to be a 
{\em $\mathbf{Q}$-factorial terminalization} of $X$ 
if $Y$ has only ${\mathbf Q}$-factorial 
terminal singularities and $f$ is a birational projective morphism 
such that $K_Y = f^*K_X$. A {\bf Q}-factorial terminalization is a 
crepant resolution exactly when $Y$ is smooth. 
Recently, Birkar-Cascini-Hacon-McKernan [B-C-H-M] have established the 
existence of minimal models of complex algebraic varieties of 
general type. As a corollary of this, we know that $X$ always has a 
{\bf Q}-factorial terminalization. In particular, $\bar{\mathcal O}$ 
should have a {\bf Q}-factorial terminalization. The author would like 
to pose the following conjecture. 
\vspace{0.15cm} 

{\bf Conjecture}. {\em Let $\mathcal{O}$ be a nilpotent orbit 
of a complex simple Lie algebra $\g$. 
Let $\tilde{\mathcal{O}}$ be the normalization of $\bar{\mathcal{O}}$. 
Then one of the following holds:} 

(1) {\em $\tilde{\mathcal{O}}$ has $\mathbf{Q}$-factorial 
terminal singularities.} 
\vspace{0.15cm}

(2) {\em There are a parabolic subalgebra $\q$ of $\g$ with  
Levi decomposition $\q = \l \oplus \n$ and a nilpotent orbit $\mathcal{O}'$ 
of $\l$ such that (a): $\mathcal{O} = \mathrm{Ind}^{\g}_{\l}(\mathcal{O}')$ 
and (b): the normalization of $G \times^Q (n(\q) + \bar{\mathcal{O}'})$ is a  
$\mathbf{Q}$-factorial terminalization of $\tilde{\mathcal{O}}$ via 
the generalized Springer map.} 

{\em Moreover, if  $\tilde{\mathcal{O}}$ does not have  $\mathbf{Q}$-factorial 
terminal singularities, then every {\bf Q}-factorial terminalization of $\tilde{\mathcal O}$ 
is of the form (2). Two {\bf Q}-factorial terminalizations are connected by Mukai flops (cf. 
[Na], p.91).}  
\vspace{0.15cm}

In this paper, we shall prove that Conjecture is true when $\g$ is 
classical.   
\vspace{0.15cm}  

{\em Acknowledgement}: The author received from B. Fu, in 2006,  a handwritten manuscript 
on explicit calculations of  the {\bf Q}-factorial terminalization for certain nilpotent orbit. 
The author would like to thank him. 
\vspace{0.15cm}

{\bf Notations and convention}.  Let $X$ be a normal algebraic variety. Then, $X$ is called {\bf Q}-{\em factorial} if, for any 
Weil divisor $D$ on $X$, its suitable multiple $mD$ ($m > 0$) is a Cartier divisor. Assume 
that the canonial divisor $K_X$ is Cartier. We say that $X$ has only {\em terminal singularities} if, for 
a resolution $\pi: Y \to X$, $K_Y = \pi^*K_X + \Sigma a_i E_i$ with $a_i > 0$ for all $i$. Here $E_i$ 
run through all $\pi$-exceptional prime divisors. For details in birational geometry, see the first part of 
[Ka] (or one can find a quick guide in [Na 4], Preliminaries).  On the other hand, for details on the basic  
properties on nilpotent orbits,  see [C-M] and [K-P]. 
        
\vspace{0.2cm}

\begin{center}
{\bf \S 1}. {\bf Preliminaries}
\end{center}

(1.1) {\em Nilpotent orbits and resolutions}:   
Let $G$ be a complex simple algebraic group and let 
$\g$ be its Lie algebra. $G$ has the adjoint 
action on $\g$. The orbit $\mathcal{O}_x$ of a nilpotent 
element $x \in \g$ for this action is called 
a nilpotent orbit.  
By the Jacobson-Morozov theorem, one can find 
a semi-simple element $h \in \g$, and a 
nilpotent element $y \in \g$ in such a way that 
$[h,x] = 2x$, $[h,y] = -2y$ and $[x,y] = h$. 
For $i \in \mathbf{Z}$, let 
$$\g_i := \{z \in \g\; [h,z] = iz \}.$$ 
Then one can write 
$$ \g = \oplus_{i \in \mathbf{Z}}\g_i. $$ 
Let $\h$ be a Cartan subalgebra of $\g$ with 
$h \in \h$.  Let $\Phi$ be the corresponding 
root system and let $\Delta$ be a base of simple roots 
such that $h$ is $\Delta$-dominant, i.e. $\alpha(h) 
\geq 0$ for all $\alpha \in \Delta$. 
In this situation, 
$$\alpha(h) \in \{0,1,2\}.$$ 
The weighted Dynkin diagram of $\mathcal{O}_x$ is 
the Dynkin diagram of $\g$ where each vertex $\alpha$ is 
labeled with $\alpha(h)$. A nilpotent orbit $\mathcal{O}_x$ 
is completely determined by its weighted Dynkin diagram. 
A Jacobson-Morozov parabolic subalgebra for 
$x$ is the parabolic subalgebra $\p$ defined by 
$$\p:= \oplus_{i \geq 0}\g_i.$$ 
Let $P$ be the parabolic subgroup of $G$ determined by  $\p$. 
We put $$ \n_2 := \oplus_{i \geq 2}\g_i. $$ 
Then $\n_2$ is an ideal of $\p$; hence,  
$P$ has the adjoint action on $\n_2$. Let us consider the   
vector bundle  $G \times^P \n_2$ over $G/P$ and the map 
$$\mu: G \times^P \n_2 \to \g$$ defined by $\mu([g,z]) := 
Ad_g(z)$. Then the image of $\mu$ coincides with the 
closure $\bar{\mathcal{O}_x}$ of $\mathcal{O}_x$ and 
$\mu$ gives a resolution of $\bar{\mathcal{O}}_x$(cf. [K-P], Proposition 7.4).   
We call $\mu$ the {\em Jacobson-Morozov resolution} 
of $\bar{\mathcal{O}_x}$. 
The orbit $\mathcal{O}_x$ has a natural closed non-degenerate 
2-form $\omega$ (cf. [C-G], Prop. 1.1.5., [C-M], 1.3). 
By $\mu$, $\omega$ is regarded as a 2-form on a Zariski open 
subset of $G \times^P \n_2$. By [Pa], [Hi], it extends to a 2-form 
on $G \times^P \n_2$. In other words, $\bar{\mathcal{O}}_x$ has 
symplectic singularity.  Let $\tilde{\mathcal O}_x$ be the normalization 
of $\bar{\mathcal O}_x$. 
In many cases, one can check the {\bf Q}-factoriality of 
$\tilde{\mathcal O}_x$ by applying the following lemma to the Jacobson-Morozov 
resolution: 
\vspace{0.15cm}

{\bf Lemma (1.1.1)}. {\em Let $\pi: Y \to X$ be a projective resolution of an 
affine variety $X$ with rational singularities. 
Let $\rho$ be the relative Picard number for $\pi$. 
If $\mathrm{Exc}(\pi)$ contains $\rho$ different prime divisors, then $X$ is 
{\bf Q}-factorial.} 

{\em Proof}. Recall that two line bundles $L$, $L'$ on $Y$ are called 
$\pi$-numerically equivalent if $(L.C) = (L'.C)$ for every proper curve $C$ 
on $Y$ such that $\pi (C)$ is a point.    
Let $N^1(\pi)$ be the group of $\pi$-numerical classes 
of line bundles on $Y$, and put  $N^1(\pi)_{\mathbf Q} := N^1(\pi) \otimes {\mathbf Q}$. 
By definition, the Picard number $\rho$ is the dimension of the {\bf Q}-vector 
space $N^1(\pi)_{\mathbf Q}$.  
Let $E_i$, ($1 \le i \le \rho$) be the prime divisors contained 
in $\mathrm{Exc}(\pi)$. We shall prove that $N^1(\pi)_{\mathbf Q} 
= \oplus \mathbf{Q}[E_i]$. 
It suffices to show that $[E_i]$'s are linearly independent. 
Put $d := \dim Y$.
Assume $\Sigma a_i [E_i] = 0$.  Put $m$ to be the largest number among $\{\dim \pi(E_i)\}$. 
We take $m$ very ample divisors $H_1$, ..., $H_m$ on $X$ in such a way that 
$\cap H_j $ intersects with $\pi (E_i)$, in finite points, for each $i$ with $\dim \pi(E_i) = m$, 
but $\cap H_j$ does not intersect with $\pi(E_i)$ for each $i$ with $\dim \pi (E_i) < m$. 
By the Bertini theorem, we may assume that $\cap \pi^*(H_j)$ is non-singular. 
We cut out $\cap \pi^*(H_j)$ further by $d - m- 2$ very ample divisors $L_1$, ..., $L_{d-m-2}$ 
on $Y$. Then the resulting variety is a smooth surface $S$. The restriction of 
$\pi$ to $S$ gives a birational contraction map $\pi\vert_S: S \to \bar{S}$. 
For $i$ with $\dim \pi(E_i) = m$, we put 
$C_i := \pi^*(H_1) \cap ... \cap \pi^*(H_m) \cap L_1\cap ... \cap L_{d-m-2} \cap E_i$.  
Then $C_i$ is a (possibly reducible) curve contained in $\mathrm{Exc}(\pi\vert_S)$.  
Moreover, two different such $C_i$'s do not have no 
common irreducible components. By a theorem of Grauert, 
the intersection matrix of  $\pi\vert_S$-exceptional curves, is 
negative definite, which implies that $[C_i]$ are linearly independent in 
$N^1(\pi\vert_S)_{\mathbf Q}$.  
Let us consider the restriction map  
$$ \iota:  N^1(\pi)_{\mathbf Q} \to N^1(\pi\vert_S)_{\mathbf Q}.$$ 
Note that $\iota ([E_i]) = 0$ for $i$ with $\dim \pi (E_i) < m$, and 
$\iota([E_i]) = [C_i]$ for $i$ with $\dim \pi (E_i) = m$.  Since 
$\Sigma a_i [E_i]) = 0$ in $N^1(\pi)_{\mathbf Q}$, we have $\Sigma a_i[C_i] = 0$ 
in $N^1(\pi\vert_S)_{\mathbf Q}$.  This implies that  
$a_i = 0$ for all $i$ with $\dim \pi (E_i) = m$.  
Next, replace $m$ by the second largest number among  $\{\dim \pi(E_i)\}$ and 
repeat the same procedure; then we finally conclude that $a_i = 0$ for all 
$i$. Now let us prove that $X$ is {\bf Q}-factorial.   
Let $D$ be a prime Weil divisor of $X$ and let $D'$ be the prime 
divisor of $Y$ obtained as the proper transform of $D$.  There are 
rational numbers $b_i$ such that $[D'] + \Sigma b_i[E_i] = 0$ in 
$N^1(\pi)_{\mathbf Q}$. Since $X$ has only rational singularities, 
$l(D' + \Sigma b_i[E_i])$ is the pull-back of a Cartier 
divisor $M$ on $X$ for some integer $l > 0$. This implies that $lD$ 
is linearly equivalent to the Cartier divisor $M$.    
\vspace{0.2cm}

(1.2) {\em Induced orbits}

(1.2.1). Let $G$ and $\g$ be the same as in (1.1). 
Let $Q$ be a parabolic subgroup of $G$ and let $\q$ be its Lie 
algebra with Levi decomposition $\q = \l \oplus \n$. Here $n$ is the 
nil-radical of $\q$ and $\l$ is a Levi-part of $\q$. 
Fix a nilpotent orbit $\mathcal{O}'$ in $\l$. Then there is a unique 
nilpotent orbit $\mathcal{O}$ in $\g$ meeting $n + \mathcal{O}'$ in 
an open dense subset ([L-S]). Such an orbit $\mathcal{O}$ is called 
the nilpotent orbit induced from $\mathcal{O}'$ and we write $$\mathcal{O} = 
\mathrm{Ind}^{\g}_{\l}(\mathcal{O}').$$  
Note that when $\mathcal{O}' = 0$, $\mathcal{O}$ is the 
Richardson orbit for $Q$. 
Since the adjoint action of 
$Q$ on $\q$ stabilizes $n + \bar{\mathcal{O}}'$, one can consider 
the variety $G \times^Q (n + \bar{\mathcal{O}}')$.   
There is a map 
$$ \nu: 
G \times^Q (n + \bar{\mathcal{O}}') \to \bar{\mathcal{O}}$$ defined 
by $\nu ([g,z]) := Ad_g(z)$. 
Since $\mathrm{Codim}_{\l}(\mathcal{O}') = 
\mathrm{Codim}_{\g}(\mathcal{O})$ (cf. [C-M], Prop. 7.1.4), 
$\nu$ is a generically finite dominating map. Moreover, $\nu$ is factorized as 
$$G \times^Q (n + \bar{\mathcal{O}}') \to G/Q \times \bar{\mathcal{O}} 
\to \bar{\mathcal{O}}$$ where the first map is a closed embedding 
and the second map is the 2-nd projection; this implies that 
$\nu$ is a projective map. In the remainder, we call $\nu$ the 
generalized Springer map for ($Q$, $\mathcal{O}'$). 

(1.2.2). Assume that $Q$ is contained in another parabolic subgroup $\bar{Q}$ 
of $G$. Let $\bar{L}$ be the Levi part of $\bar{Q}$ which contains 
the Levi part $L$ of $Q$. Let $\bar{\q} = \bar{\l}\oplus\bar{\n}$ be the 
Levi decomposition.  
Note that $\bar{L} \cap Q$ is a parabolic subgroup of 
$\bar{L}$ and $\l(\bar{L} \cap Q) = \l$. Let  
$\mathcal{O}_1 \subset \bar{\l}$ be the nilpotent orbit induced from $(\bar{L} \cap Q, \mathcal{O}')$. 
Then there is a natural map 
$$\pi: G \times^Q (n + \bar{\mathcal{O}}') \to G \times^{\bar{Q}} (\bar{n} + \bar{\mathcal{O}_1})$$ 
which factorizes $\nu$ as $\bar{\nu}\circ\pi = \nu$. Here $\bar{\nu}$ is the generalized 
Springer map for $(\bar{Q}, \mathcal{O}_1)$. 

(1.2.3). Assume that there are a parabolic subgroup $Q_L$ of $L$ and 
a nilpotent orbit $\mathcal{O}_2$ in the Levi subalgebra $\l(Q_L)$ such that 
$\mathcal{O}'$ is the nilpotent orbit induced from $(Q_L, \mathcal{O}_2)$. 
Then there is a parabolic subgroup $Q'$ of $G$ such that $Q' \subset Q$, 
$\l(Q') = \l(Q_L)$ and $\mathcal{O}$ is the nilpotent orbit induced from 
$(Q', \mathcal{O}_2)$. The generalized Springer map $\nu'$ for $(Q', \mathcal{O}_2)$ 
is factorized as 
$$ G \times^{Q'}(\n' + \bar{\mathcal{O}}_2) \to G \times^Q(\n + \bar{\mathcal{O}'}) \to \bar{\mathcal{O}}.$$ 

{\bf  Lemma (1.2.4)}. {\em Let 
$$ \nu: G \times^Q (n + \bar{\mathcal{O}}') \to \bar{\mathcal{O}}$$ be a generalized Springer 
map defined in (1.2.1). Then the normalization of $G \times^Q (n + \bar{\mathcal{O}}')$ is a symplectic variety.} 
\vspace{0.15cm}

{\em Proof}.  We shall prove that the pull-back of the Kostant-Kirillov form $\omega$ on $\mathcal{O}$ 
gives a non-degenerate 2-form on $G \times^Q (n + {\mathcal{O}}')$. 
This is enough for proving (1.2.4).  
In fact, $G \times^Q (n + \bar{\mathcal{O}}')$ is locally a product of 
$\bar{\mathcal{O}}'$ and a non-singular variety; hence 
the normalization of $G \times^Q (n + \bar{\mathcal{O}}')$ has only rational 
Gorenstein singularities. The Kostant-Kirillov form extends to a regular 2-form on  
any resolution of $G \times^Q (n + \bar{\mathcal{O}}')$ as explained in (1.1). 
Therefore, the normalization of $G \times^Q (n + \bar{\mathcal{O}}')$ is a symplectic variety. 
The following argument is analogous to [Pa]. Let $\h$ be a Cartan subalgebra of $\g$ such that 
$\h \subset \l$. There is an involution $\phi_{\g}$ of $\g$ which stabilizes 
$\h$ and which acts on the root system $\Phi$ via $-1$. Put $\n_{-} := \phi_{\g}(\n)$. 
Take a point $[1, y+y'] \in G \times^Q (n + {\mathcal{O}}')$ so that $y \in \n$, $y' \in \mathcal{O}'$ 
and $y + y' \in \mathcal{O}$. The tangent space of 
$G \times^Q (n + \bar{\mathcal{O}}')$ at $[1, y+y']$ is decomposed as 
$$ T_{[1, y+y']} = \n_{-} \oplus T_{y+y'}(\n + \bar{\mathcal O}'). $$ 
Since $Q\cdot (y+y')$ coincides with the Zariski open dense subset 
$\mathcal{O} \cap (\n + {\mathcal O}') \subset \n + {\mathcal O}'$,  
an element $v \in T_{[1, y+y']}$ can be written as 
$$ v = v_1 + [v_2, y+y'], \;\; v_1 \in \n_{-}, \; v_2 \in \q.$$ 
Let $d\nu_*: T_{[1, y+y']} \to T_{\nu([1, y+y'])}\mathcal{O}$ be the tangential map 
for $\nu$. 
Then $$d\nu_*(v) = [v_1 + v_2, y+y'].$$ 
Take one more element $w \in T_{[1, y+y']}$ in such a way that 
$$ w = w_1 + [w_2, y+y'], \;\; w_1 \in \n_{-}, \; w_2 \in \q.$$ 
Denote by $\langle \;, \;\rangle$ the Killing form of $\g$. 
By the definition of the Kostant-Kirillov form, one has 
$$ \omega(d\nu_*(v), d\nu_*(w)) := \langle y+y', [v_1 + v_2, w_1 + w_2]\rangle.$$ 
Note that $\langle y+y', [v_1, w_1]\rangle = \langle y, [v_1, w_1]\rangle$, and 
$\langle y+y', [v_2, w_2]\rangle = \langle y', [v_2, w_2]\rangle$.  Therefore, 
\vspace{0.2cm}

$\omega(d\nu_*(v), d\nu_*(w)) =  $ 

$\langle y, [v_1,w_1] \rangle +  \langle y+y', [v_1, w_2] \rangle   
+ \langle y+y', [v_2, w_1] \rangle + \langle y', [v_2, w_2] \rangle = $ 

$\langle y, [v_1,w_1] \rangle + \langle v_1, [w_2, y + y']_n\rangle  - 
\langle [v_2, y + y']_n, w_1 \rangle  
+ \omega ([v_2, y'], [w_2, y']), $
\vspace{0.2cm}

where 
$[w_2, y + y']_n$ (resp. $[v_2, y + y']_n$) is 
the nil-radical part of  $[w_2, y + y']$ (resp. 
$[v_2, y + y']$) in the decomposition 
$T_{y+y'}(\n + \bar{\mathcal O}') = \n + T_{y'}\mathcal{O}'$. 
Let $\mathcal{O}_r \subset \g$ be the Richardson orbit for $Q$, and 
let $\pi: T^*(G/Q) \to \bar{\mathcal{O}_r}$ be the Springer map. 
The first part $\langle y, [v_1,w_1] \rangle + \langle v_1, [w_2, y + y']_n\rangle  - 
\langle [v_2, y + y']_n, w_1 \rangle$ 
corresponds to the 2-form on $T^*(G/Q)$ obtained by the pull-back of 
the Kostant-Kirillov 2-form on $\mathcal{O}_r$ by $\pi$ (cf. [Pa]), which is non-degenerate 
on $T^*(G/Q)$. 
Let us consider the second part $\omega ([v_2, y'], [w_2, y'])$. 
Denote by $[v_2, y+y']_l$ (resp. 
$[w_2, y+y']_l$) the $T_{y'}\mathcal{O}'$-part of $[v_2, y+y']$ (resp. $[w_2, y+y']$) in 
the decomposition 
$T_{y+y'}(\n + \bar{\mathcal O}') = \n + T_{y'}\mathcal{O}'$. 
Then 
$[v_2,y'] = [v_2, y+y']_l$ and $[w_2, y'] = [w_2, y+y']_l$; hence,  
the second part is the Kostant-Kirillov form on $\mathcal{O}'$.  
The arguments above show that      
$\nu^*\omega$ is non-degenerate at $[1, y+y']$ for an arbitary 
$y \in \n$ and for an arbitrary $y' \in {\mathcal O}'$. 
By the $G$-equivariance of $\nu$, we have the lemma. 
\vspace{0.15cm}

(1.3) {\em Nilpotent orbits in classical Lie algebras}: 
When $\g$ is a classical Lie algebra, $\g$ is naturally 
a Lie subalgebra of $\mathrm{End}(V)$ for a complex vector space 
$V$. Then we can attach a partition ${\bf d}$ of 
$n := \dim V$ to each orbit as the Jordan 
type of an element contained in the orbit. 
Here a 
partition ${\bf d} := [d_1, d_2,..., d_k]$ of $n$ is a set of 
positive integers with $\Sigma d_i = n$ and 
$d_1 \geq d_2  \geq ... \geq d_k$. 
Another way of writing $\mathbf{d}$ is 
$[d_1^{s_1}, ..., d_k^{s_k}]$ with $d_1 > d_2 ... > d_k >0$. 
Here $d_i^{s_i}$ is an $s_i$ times $d_i$'s: 
$d_i, d_i, ..., d_i$.   
The partition  ${\mathbf d}$ corresponds to a Young diagram. 
For example, $[5,  4^2,  1]$ corresponds to 

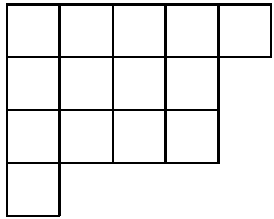
\begin{picture}(100, 100)(0, 0)
\put(00,  80){\line(1, 0){100}}
\put(00,  60){\line(1, 0){100}}
\put(00,  40){\line(1, 0){80}}
\put(00,  20){\line(1, 0){80}}
\put(00,  00){\line(1, 0){20}}

\put(00,  00){\line(0, 1){80}}
\put(20,  00){\line(0, 1){80}}
\put(40,  20){\line(0, 1){60}}
\put(60,  20){\line(0, 1){60}}
\put(80,  20){\line(0, 1){60}}
\put(100,  60){\line(0, 1){20}}
\end{picture}

When an   
integer $e$ appears in the partition $\bf{d}$, we say 
that $e$ is a {\em member} of $\bf{d}$. We call $\bf{d}$ 
{\em very even} when $\bf{d}$ consists with only even 
members, each having even multiplicity. 

Let us denote by $\epsilon$ the number $1$ or $-1$. 
Then a partition $\mathbf{d}$ is $\epsilon$-admissible 
if all even (resp. odd) members of $\mathbf{d}$ have  
even multiplicities when $\epsilon = 1$ (resp. $\epsilon = 
-1$).  
The following result can be found, for example, in 
[C-M, \S 5].     

{\bf Proposition (1.3.1)}
{\em   Let $\mathcal{N}o(\g)$ be the set of nilpotent 
orbits of $\g$.} 
\vspace{0.12cm}

(1)($A_{n-1}$): {\em When $\g = \mathfrak{sl}(n)$, there is a 
bijection between $\mathcal{N}o(\g)$ and 
the set of partitions $\bf{d}$ of $n$. }  
\vspace{0.12cm} 
 
(2)($B_n$): {\em When $\g = \mathfrak{so}(2n+1)$, there is a 
bijection between $\mathcal{N}o(\g)$ and 
the set of $\epsilon$-admissible partitions $\bf{d}$ of $2n+1$ 
with $\epsilon = 1$.}    
\vspace{0.12cm}

(3)($C_n$): {\em When $\g = \mathfrak{sp}(2n)$, there is a 
bijection between $\mathcal{N}o(\g)$ and the 
set of $\epsilon$-admissible partitions $\bf{d}$ of $2n$ 
with $\epsilon = -1$.}
\vspace{0.12cm}

(4)($D_n$): {\em When $\g = \mathfrak{so}(2n)$, there is a 
surjection $f$ from $\mathcal{N}o(\g)$ to the set 
of $\epsilon$-admissible partitions $\bf{d}$ of $2n$ with 
$\epsilon = 1$. For a partition $\bf{d}$ which is 
not very even, $f^{-1}(\bf{d})$ consists of exactly one orbit, 
but, for very even $\bf{d}$, $f^{-1}(\bf{d})$ consists of exactly 
two different orbits. }  \vspace{0.2cm}

Take an $\epsilon$-admissible partition $\mathbf{d}$ of a 
positive integer $m$. If $\epsilon = 1$, we put $\g = so(m)$ 
and if $\epsilon = -1$, we put $\g = sp(m)$.   
We denote by $\mathcal{O}_{\mathbf d}$ a nilpotent orbit in $\g$ 
with Jordan type $\mathbf{d}$. Note that, except when $\epsilon = 1$  
and $\mathbf{d}$ is very even, $\mathcal{O}_{\mathbf{d}}$ is uniquely 
determined. When $\epsilon = 1$ and $\mathbf{d}$ is very even, 
there are two possibilities for $\mathcal{O}_{\mathbf{d}}$. If necessary, 
we distinguish the two orbits by the labelling: $\mathcal{O}^I_{\mathbf d}$ 
and $\mathcal{O}^{II}_{\mathbf d}$. 
Let us fix a classical Lie algebra $\g$ and study the relationship 
among nilpotent orbits in $\g$. 
When $\g$ is of type $B$ or $D$ (resp. $C$), we only consider 
the $\epsilon$-admissible partitions with $\epsilon = 1$ (resp. 
$\epsilon = -1$). 
We introduce a partial order in the set of the  
partitions of (the same number):  
for two partitions $\mathbf{d}$ and $\mathbf{f}$,  
$\mathbf{d} \geq \mathbf{f}$ if $\Sigma_{i \leq k}d_i 
\geq \Sigma_{i \leq k}f_i$ for all $k \geq 1$.  
On the other hand, for two nilpotent orbits $\mathcal{O}$ 
and $\mathcal{O}'$ in  $\g$, we write $\mathcal{O} \geq {\mathcal O}'$ 
if $\mathcal{O}' \subset \bar{\mathcal{O}}$. 
Then, $\mathcal{O}_{\mathbf d} \geq \mathcal{O}_{\mathbf f}$ if and 
only if $\mathbf{d} \geq \mathbf{f}$. When $\mathbf{d}$ and 
$\mathbf{f}$ are $\epsilon$-admissible partitions with 
$\mathbf{f} \geq \mathbf{g}$, we call this pair an {\em $\epsilon$-degeneration} 
or simply a {\em degeneration}. 

Now let us consider the case $\g$ is of type $B$, $C$ or $D$. 

Assume that an $\epsilon$- degeneration $\mathbf{d} \geq \mathbf{f}$ is  
{\em minimal} in the sense that  
there is no $\epsilon$-admissible partition 
$\mathbf{d}'$ (except $\mathbf{d}$ and $\mathbf{f}$) 
such that $\mathbf{d} \geq \mathbf{d}' \geq \mathbf{f}$. 
Kraft and Procesi [K-P] have studied the 
{\em normal slice} $N_{\mathbf{d}, \mathbf{f}}$ of 
$\mathcal {O}_{\mathbf f} \subset \bar{\mathcal{O}}_{\mathbf d}$ 
in such cases. 
 If, for two integers $r$ and $s$, the first $r$ rows and the first 
$s$ columns of $\mathbf{d}$ and $\mathbf{f}$ coincide and the 
partition $(d_1, ..., d_r)$ is  $\epsilon$-admissible, 
then one can erase these rows and columns from $\mathbf{d}$ and 
$\mathbf{f}$ respectively to get new partitions $\mathbf{d}'$ 
and $\mathbf{f}'$ with $\mathbf{d}' \geq \mathbf{f}'$. 
If we put $\epsilon' := (-1)^s\epsilon$, then $\mathbf{d}'$ 
and $\mathbf{f}'$ are both $\epsilon'$-admissible. 
The pair $(\mathbf{d}', \mathbf{f}')$ is also minimal. 
Repeating such process, one can reach a degeneration 
$\mathbf{d}_{irr} \geq \mathbf{f}_{irr}$ which is {\em irreducible} 
in the sense that there are no rows and columns to be erased.  
By [K-P], Theorem 2, $N_{\mathbf{d}, \mathbf{f}}$ is analytically 
isomorphic to $N_{\mathbf{d}_{irr}, \mathbf{f}_{irr}}$ around the 
origin. 
According to [K-P], a minimal and irreducible degeneration $\mathbf{d} \geq  \mathbf{f}$ 
is one of the following: 

a: $\g = sp(2)$, $\mathbf{d} = (2)$, and $\mathbf{f} = (1^2)$. 

b: $\g = sp(2n)$ ($n > 1$), $\mathbf{d} = (2n)$, and $\mathbf{f} = (2n-2, 2)$.

c: $\g = so(2n+1)$ ($n > 0$), $\mathbf{d} = (2n+1)$, and $\mathbf{f} = (2n-1, 1^2)$. 

d: $\g = sp(4n+2)$ ($n > 0$), $\mathbf{d} = (2n+1, 2n+1)$, and $\mathbf{f} = (2n,2n,2)$.

e: $\g = so(4n)$ ($n > 0$), $\mathbf{d} = (2n, 2n)$, and $\mathbf{f} = (2n-1,2n-1, 1^2)$.

f: $\g = so(2n+1)$ ($n > 1$), $\mathbf{d} = (2^2, 1^{2n-3})$, and $\mathbf{f} = (1^{2n+1})$.

g: $\g = sp(2n)$ ($n > 1$), $\mathbf{d} = (2, 1^{2n-2})$, and $\mathbf{f} = (1^{2n})$. 

h: $\g = so(2n)$ ($n > 2$), $\mathbf{d} = (2^2, 1^{2n-4})$, and $\mathbf{f} = (1^{2n})$.

In the first 4 cases (a,b,c,d,e),  $\mathcal{O}_{\mathbf{f}}$ have codimension 2 
in $\bar{\mathcal{O}}_{\mathbf{d}}$. 
In the last 3 cases (f,g,h), $\mathcal{O}_{\mathbf{f}}$ have codimension $ \geq 4$  
in $\bar{\mathcal{O}}_{\mathbf{d}}$. 
\vspace{0.2cm}

{\bf  Proposition (1.3.2)} 
{\em Let $\mathcal{O}$ be a nilpotent orbit in a classical Lie algebra $\g$ 
of type $B$, $C$ or $D$ with Jordan type $\mathbf{d} := [(d_1)^{s_1}, ..., (d_k)^{s_k}]$ 
($d_1 >d_2 > ... > d_k$).  Let $\Sigma$ be 
the singular locus of  $\bar{\mathcal{O}}$. 
Then $\mathrm{Codim}_{\bar{\mathcal{O}}}(\Sigma) \geq 4$ if and only if 
the partition $\mathbf{d}$ has full members, that is,  
any integer $j$ with $1 \leq j \leq d_1$ is a member of $\mathbf{d}$. 
Otherwise, $\mathrm{Codim}_{\bar{\mathcal{O}}}(\Sigma) = 2$. } 
\vspace{0.2cm}

{\em Proof}. 
Assume that $\mathrm{Codim}_{\bar{\mathcal{O}}}(\Sigma) \geq 4$. 
We shall prove that $\mathbf{d}$ has full members. Suppose, to the 
contrary, that there is some $i$ with $d_i \geq d_{i+1} +2$.  
Then one can find a minimal degeneration $\mathbf{d} \geq \mathbf{f}$ 
where $\mathbf{f}$ is one of the following:
$$[..., d_i^{s_i-1}, (d_{i+1}+1)^2, d_{i+1}^{s_{i+1}-1}, ...], \;\; d_i = d_{i+1} + 2$$ 
$$[..., d_i^{s_i-1}, d_i-2, d_{{i+1}}+2, d_{{i+1}}^{s_{i+1}-1}, ...]$$
$$[..., d_i^{s_i-1}, d_i-2, (d_{{i+1}}+1)^2, d_{{i+1}}^{s_{i+1}-2}, ...]$$ 
$$[..., d_i^{s_i-2}, (d_i-1)^2, d_{{i+1}}+2, d_{{i+1}}^{s_{i+1}-1}, ...]$$ 
$$[..., d_i^{s_i-2}, (d_{i}-1)^2, (d_{{i+1}}+1)^2, d_{{i+1}}^{s_{i+1}-2}, ...].$$ 
By the row-column erasing one gets  
an irreducible, minimal degeneration $\mathbf{d}_{irr} \geq \mathbf{f}_{irr}$ 
of type $a$,$b$,$c$,$d$ or $e$. This is a contradiction; hence 
$\mathbf{d}$ has full members.  
Conversely, assume that $\mathbf{d}$ has full members, i.e. $d_j = k-j+1$ for all 
$1 \leq j \leq k$. Then, for every minimal degeneration   
$\mathbf{d} \geq \mathbf{f}$, $\mathbf{f}$  
is one of the following:  
$$[..., (k-i+1)^{s_i-2}, (k-i)^{s_{i-1}+4}, (k-i-1)^{s_{i-2}-2}, ...]$$ 
$$[..., (k-i+1)^{s_i-1}, (k-i)^{s_{i-1}+2}, (k-i-1)^{s_{i-2}-1}, ...].$$
By the row-column erasing one gets  
an irreducible, minimal degeneration $\mathbf{d}_{irr} \geq \mathbf{f}_{irr}$ 
of type $f$,$g$ or $h$. Therefore,  
$\mathrm{Codim}_{\bar{\mathcal{O}}}(\Sigma) \geq 4$.    
\vspace{0.2cm}

(1.4) {\em Jacobson-Morozov resolutions in the case of classical Lie algebras}(cf. 
[CM], 5.3): 
Let $V$ be a complex vector space of dimension $m$ with a non-degenerate 
symmetric (or skew-symmetric) form $<\;, \; >$. In the symmetric case,  we take a 
basis $\{e_i\}_{1 \le i \le m}$ of $V$ in such a way that $<e_j, e_k> = 1$ if $j + k = m+1$ 
and otherwise $<e_j, e_k> = 0$. In the skew-symmetric case, we take a basis 
$\{e_i\}_{1 \le i \le m}$ of $V$ in such a way that $<e_j, e_k> = 1$ if $j < k$ and $j + k = m+1$, 
and $<e_j, e_k> = 0$ if $j + k \ne m+1$.  
When $(V, <\;, \;>)$ is a symmetric vector space, $\g := so(V)$ is the Lie algebra 
of type $B_{(m-1)/2}$ (resp. $D_{m/2}$)  if $m$ is odd (resp. even). 
When $(V, <\;, \;>)$ is a skew-symmetric vector space, $\g := sp(V)$ is the Lie algebra 
of type $C_{m/2}$. In the remainder of this paragraph, $\g$ is one of these Lie 
algebra contained in $\mathrm{End}(V)$. 
Let $\h \subset \g$ be the Cartan subalgebra consisting of all diagonal 
matrices, and let $\Delta$ be the standard base of simple roots.  
Let $x \in \g$ be a nilpotent element. As in (1.1), one can choose 
$h$, $y \in \g$ in such a way that $\{x,y,h\}$ is a $sl(2)$-triple. 
If necessary, by replacing $x$ by its conjugate element, one may assume 
that $h \in \h$ and $h$ is $\Delta$-dominant.  
Assume that $x$ has Jordan type 
$\mathbf{d} = [d_1, ...., d_k]$.   
The diagonal matrix $h$ is described as follows. 
Let us consider the sequence of integers of length $m$: \vspace{0.2cm}
 
$d_1-1, d_1-3 , ..., -d_1 + 3, -d_1 +1, d_2-1, d_2-3, ..., -d_2 +3, -d_2 + 1, 
..., d_k -1, d_k -3, ..., -d_k +3, -d_k +1$.   
\vspace{0.2cm}
       
Rearrange this sequence in the non-increasing order and get a new sequence 
$p_1^{t_1}, ..., p_l^{t_l}$ with $p_1 > p_2 ... > p_l$ and $\Sigma t_i = m$.    
Then $$h = \mathrm{diag}(p_1^{t_1}, ..., p_l^{t_l}).$$ 
Here $p_i^{t_i}$ means the $t_i$ times of $p_i$'s: $p_i, p_i, ..., p_i$. 
It is then easy to describe explicitly the Jacobson-Morozov parabolic subalgebra 
$\p$ of $x$ and its ideal $\n_2$ (cf. (1.1)). 
The Jacobson-Morozov parabolic subgroup $P$ is the stabilizer group of 
certain isotropic flag $\{F_i\}_{1 \le i \le r}$ of $V$. 
Here, an isotropic flag of $V$ (of length $r$) is a increasing filtration 
$0 \subset F_1 \subset F_2 \subset ... \subset F_r \subset V$ such that 
$F_{r+1-i} = F_i^{\perp}$ for all $i$. The flag type of $P$ is $(t_1, ..., t_l)$. 
The nilradical $\n := \oplus_{i > 0}\g_i$ of $\p$ consists of the 
elements $z$ of $\g$ such that $z(F_i) \subset F_{i-1}$ for all $i$. 
On the other hand, it depends on the weighted Dynkin diagram for $x$ 
how $\n_2$ takes its place in $\n$. 
\vspace{0.2cm}

{\bf Example (1.4.1)}.   
Assume that $\mathbf{d}$ has full members (cf. (1.3.2)). 
Then, it can be checked that 
$$ \n_2 = \{z \in \g ; z(F_i) \subset F_{i-2}\;\; \mathrm{for} \;\; \mathrm{all} \;  i\}. $$ 
Let us consider   
the Jacobson-Morozov resolution  
$$ \mu: G \times^P \n_2 \to \bar{\mathcal{O}}_{\mathbf{d}}.$$ 
Take $z \in \mathcal{O}_{\mathbf{d}}$. One can look at the fiber $\mu^{-1}(z)$ 
by using the characterization of $\n_2$ above. 
Assume that $G = Sp(V)$. We prepare two kinds of skew-symmetric 
vector space $V_d$ ($d$: even), and $W_{2d}$ ($d$: odd) as follows. 
The $V_d$ is a $d$ dimensional vector space with 
a basis $e_{d-1}$, $e_{d-3}$, ..., $e_{-d+3}$, $e_{-d+1}$. 
The skew-symmetric form $\langle \;,\; \rangle$ is defined in such a way that 
$\langle e_i, e_j \rangle = 1$
(resp. $\langle e_i,e_j\rangle = -1$)  when $i + j = 0$ and $i >j$ (resp. $i+j = 0$ and $i < j$), and 
$\langle e_i, e_j\rangle = 0$ when $i + j \ne 0$. 
Let $z_d$ be the endomorphism of $V_d$ defined by  the $d \times d$ matrix 
$J_d$ with $J_d(i,i+1) = 1$ ($1 \leq i \leq d/2$), $J_d(i,i+1) = -1$ ($d/2 + 1 \leq i \leq d-1$) 
and otherwise $J_d(i,j) = 0$.   
The $W_{2d}$ is a $2d$ dimensional vector space with a basis 
$f_{d-1}$, $f_{d-3}$, ..., $f_{-d+3}$, $f_{-d+1}$, 
$f'_{d-1}$, $f'_{d-3}$, ..., $f'_{-d+3}$, $f'_{-d+1}$. 
The skew-symmetric form $\langle \;,\;\rangle$ is defined in such a way that 
$\langle f_i, f'_j\rangle = 1$
($\langle f'_i, f_j\rangle = -1$) when $i + j = 0$, and otherwise $\langle f_i, f'_j\rangle = 0$, 
$\langle f'_i, f_j\rangle = 0$.  
Let $z_{2d}$ be the endomorphism defined by the matrix 
$J'_{2d}$ with $J'_{2d}(i, i+1) = 1$ ($1 \le i \le d-1$), $J'_{2d}(d,d+1) = 0$, 
$J'_{2d}(i,i+1) = -1$ ($d+1 \le i \le 2d-1$) and otherwise $J'_{2d}(i,j) = 0$.   
Write  $\mathbf{d}$ as  $[k^{s_k}, (k-1)^{s_{k-1}}, ..., 2^{s_2}, 1^{s_1}]$ 
with $s_i > 0$. Note that when $i$ is odd, $s_i$ is even. 
In the notation above, $d_1 = ... = d_{s_k} = k$, 
$d_{s_k+1} = ... = d_{s_k + s_{k-1}} = k-1$ and so on. 
Then, $t_1 = s_k$, $t_2 = s_{k-1}$, $t_3 = s_k + s_{k-2}$, 
$t_4 = s_{k-1} + s_{k-3}$, ...
We may assume that $$V = (\bigoplus_{d: even} (V_d)^{\oplus s_d}) \oplus  
(\bigoplus_{d: odd}(W_{2d})^{\oplus s_d/2})$$ and 
$z\vert_{V_d} = z_d$, $z\vert_{W_{2d}} = z_{2d}$. 
An element of $\mu^{-1}(z)$ is the isotropic flag $\{F_i\}$ of $V$ 
with flag type $(t_1, ..., t_l)$ which satisfies 
$z(F_i) \subset F_{i-2}$ for all $i$.   
One can find such a flag as follows.   
Put $F_1 := z^{k-1}(V) (= \mathrm{Im}(z^{k-1}))$ and define 
$F_{2k-2} := F_1^{\perp}$. We next put 
$F_2 := z^{k-2}(F_{2k-2})$ and $F_{2k-3} := F_2^{\perp}$. 
The subsequent step is similar;  we put inductively $F_i := z^{k-i}(F_{2k-i})$ 
and $F_{2k-i-1} := F_i^{\perp}$.  

When $G = SO(V)$, $\mathbf{d}$ can be written as  $[k^{s_k}, (k-1)^{s_{k-1}}, ..., 2^{s_2}, 1^{s_1}]$ 
where $s_i$ is even when $i$ is even. We prepare two kinds of 
symmetric vector spaces $V_d$ ($d$: odd) and $W_{2d}$ ($d$: even). 
Then $$V = (\bigoplus_{d: odd} (V_d)^{\oplus s_d}) \oplus  
(\bigoplus_{d: even}(W_{2d})^{\oplus s_d/2}).$$ 
The description of the flag corresponding to $\mu^{-1}(z)$ is 
quite similar. 
\vspace{0,2cm}

{\bf Lemma (1.4.2)} {\em Assume that $\mathbf{d}$ has full members. 
For each minimal $\epsilon$-degeneration $\mathbf{d} \geq \mathbf{f}$, the fiber 
$\mu^{-1}(\mathcal{O}_{\mathbf{f}})$ has codimension $1$ in 
 $G \times^P \n_2$.} 
\vspace{0.2cm}

{\em Proof}. {\bf (1)}: By the proof of (1.3.2), if $\mathbf{d}$ has full members, 
then, for every minimal degeneration $\mathbf{d} \geq \mathbf{f}$, its 
reduction $\mathbf{d}_{irr} \geq \mathbf{f}_{irr}$ is of type $f$, $g$ 
or $h$. 
If it is of type $f$, the normal slice $N_{\mathbf{d}, \mathbf{f}}$ of 
$\mathcal{O}_{\mathbf{f}} \subset \bar{\mathcal{O}}_{\mathbf{d}}$ is  
isomorphic to the germ of $\bar{\mathcal{O}}_{[2^2, 1^{2n-3}]} (\subset so(2n+1))$ 
($n > 1$) at $0$. 
Similarly, if it is of type $g$ (resp. $h$),  
$N_{\mathbf{d}, \mathbf{f}}$ is isomorphic to the germ of 
$\bar{\mathcal{O}}_{[2, 1^{2n-2}]} (\subset sp(2n))$ ($n > 1$)  
(resp. $\bar{\mathcal{O}}_{[2^2, 1^{2n-4}]} (\subset so(2n))$ ($n > 2$)) 
at $0$. 
Note that they are all isolated singularities. 
Except when $\mathbf{d}_{irr} \geq \mathbf{f}_{irr}$ is of type 
$h$ with $n = 3$, these germs have $\mathbf{Q}$-factorial terminal 
singularities. Indeed, they are isolated symplectic singularities of dim $\geq 4$; 
hence they have only terminal singularities. The $\mathbf{Q}$-factoriality 
of them is checked by using the Jacobson-Morozov resolutions of them. 
The exceptional locus of each Jacobson-Morozov resolution consists 
of a flag variety with $b_2 = 1$; this implies the $\mathbf{Q}$-factoriality. 
If $\mathrm{Codim}\;\mu^{-1}({\mathcal{O}}_{\mathbf{f}}) 
\geq 2$, then $N_{\mathbf{d}, \mathbf{f}}$ is not $\mathbf{Q}$-factorial. 
Therefore, the lemma follows in these cases. 

{\bf (2)}: 
The only exception is when $\mathbf{d}_{irr} \geq \mathbf{f}_{irr}$ is of type 
$h$ with $n = 3$.   
In this case, the germ of $\bar{\mathcal O}_{[2^2, 1^2]} (\subset so(6))$ at 
$0$ has only terminal singularity because it is an isolated symplectic singularity 
with $\dim 6$.  But its Jacobson-Morozov resolution has $\mathrm{Gr}_{iso}(2,6)$ 
as its exceptional divisor and $b_2(\mathrm{Gr}_{iso}(2,6)) = 2$. 
By this observation, we know that $\bar{\mathcal O}_{[2^2, 1^2]}$ is 
not $\mathbf{Q}$-factorial.  
Let $\mathbf{d} = [k^{s_k}, (k-1)^{s_{k-1}}, ..., 2^{s_2}, 1^{s_1}]$ be an 
$\epsilon$-admissible partition of $m$.   
Our exceptional case only occurs when $s_j = 2$ for some $j$ 
with $(-1)^{j+1} = \epsilon$.  Then 
$$ \mathbf{f} = [..., (j+1)^{s_{j+1}-2}, j^6, (j-1)^{s_{j-1}-2}, ...].$$ 
Let  $$ \mu: G \times^P \n_2 \to \bar{\mathcal{O}}_{\mathbf{d}}$$
be the Jacobson-Morozov resolution.
Take $z \in \mathcal{O}_{\mathbf{f}}$. We shall prove that 
the isotropic Grassmann variety $Gr_{iso}(2,6)$ is contained in 
$\mu^{-1}(z)$ by a similar argument to (1.4.1). 
We only discuss here the case when $G = Sp(V)$, but the case 
of $G = SO(V)$ is quite similar. 
As in (1.4.1), we take two skew-symmetric spaces $V_d$ ($d$: even) 
and $W_{2d}$ ($d$: odd).    
We may assume that $V$ is isomorphic to $$\bigoplus_{d: even, d \ne j-1, j, j+1} (V_d)^{\oplus s_d} \oplus 
W_{2j-2}^{\oplus s_{j-1}/2-1} \oplus V_j^{\oplus 6} \oplus W_{2j+2}^{\oplus s_{j+1}/2-1} \oplus   
\bigoplus_{d: odd, d \ne j-1, j, j+1}(W_{2d})^{\oplus s_d/2}$$ and 
$z\vert_{V_d} = z_d$, $z\vert_{W_{2d}} = z_{2d}$. 
We  must find isotropic flags $\{F_i\}$ of $V$ 
with flag type $(t_1, ..., t_l)$ which satisfies 
$z(F_i) \subset F_{i-2}$ for all $i$. Here $(t_1, ..., t_l)$ is 
the same one as in (1.4.1).  As in (1.4.1), we define $F_i$ 
with $i \leq k-j-1$ by  
$F_i := z^{k-i}(F_{2k-i})$ 
and $F_{2k-i-1} := F_i^{\perp}$. 
But, the situation is different from (1.4.1) when $i = k-j$. 
We cannot put $F_{k-j} := z^j(F_{k+j})$. 
In fact, $\dim z^j(F_{k+j})$ is exactly 2 less than the dimension 
$F_{k-j}$ should have because the exponent $s_{j+1} -2$ 
of $j+1$ in $\mathbf{f}$ 
is exactly 2 less than that of $j+1$ in $\mathbf{d}$. 
Let us consider the direct summand 
$V_j^{\oplus 6}$ of $V$.  The kernel of the endomorphism 
$z_j: V_j \to V_j$ has dimension 1 and is spanned by $e_{j-1}$. 
Take 6 copies $e_{j-1}^{(1)}, ..., e_{j-1}^{(6)}$ of $e_{j-1}$.
Then, $\mathrm{Ker}(z_j^{\oplus 6})$ is a 6 dimensional 
vector space spanned by  $e_{j-1}^{(1)}, ..., e_{j-1}^{(6)}$. 
We want to choose two dimensional subspace $L \subset 
\mathrm{Ker}(z_j^{\oplus 6})$ and want to define 
$F_{k-j}$ as $z^j(F_{k+j}) + L$. Since $F_{k-j}$ should 
be contained in $z^{j-1}((F_{k-j})^{\perp})$,  
we must choose $L$ in such a way that 
$$z^j(F_{k+j}) + L \subset z^{j-1}((z^j(F_{k+j}) + L)^{\perp}).$$ 
Let $v = \Sigma a_i e_{j-1}^{(i)}$ and $w = \Sigma b_i e_{j-1}^{(i)}$ 
be a basis of $L$. Then the condition above is equivalent 
to 
$$\Sigma a_i^2 = \Sigma b_i^2 = \Sigma a_ib_i = 0.$$ 
This means that $[L] \in Gr_{iso}(2,6)$. 
For such an $L$, we put $F_{k-j} := z^j(F_{k+j}) + L$. 
Once $F_{k-j}$ is fixed, we define, for $i \geq k-j+1$, 
$F_i := z^{k-i}(F_{2k-i})$ 
and $F_{2k-i-1} := F_i^{\perp}$.
One can check that $\{F_i\}$ is a desired flag. 
Therefore, $Gr_{iso}(2,6) \subset  
\mu^{-1}(z)$. 
Since $\dim Gr_{iso}(2,6) = 5$ and 
$\dim N_{\mathbf{d}, \mathbf{f}} = 6$, this implies that    
$\mu^{-1}({\mathcal{O}}_{\mathbf{f}})$ has codimension $1$ in  
 $G \times^P \n_2$.  
\vspace{0.2cm}

{\bf Corollary (1.4.3)} {\em Assume that $\mathbf{d}$ is an $\epsilon$-admissible 
partition and it has full members. Let $\tilde{\mathcal{O}}_{\mathbf{d}}$ be 
the normalization of $\bar{\mathcal{O}}_{\mathbf{d}}$. Then, 
$\tilde{\mathcal{O}}_{\mathbf{d}}$ has only $\mathbf{Q}$-factorial termainal 
singularities except when $\g = so(4n+2)$, $n \geq 1$ and $\mathbf{d} = [2^{2n}, 1^2]$.}  

{\em Proof}. 
Let $k$ be the maximal member of $\mathbf{d}$. 
Then there are $k-1$ minimal degenerations $\mathbf{d} \geq \mathbf{f}$. 
By Lemma (1.4.2), $\mathrm{Exc}(\mu)$ 
contains at least $k-1$ irreducible divisors. When $\epsilon = 1$ (i.e, 
$\g = so(V)$) and 
there is a minimal degeneration $\mathbf{d} \geq \mathbf{f}$ with 
$\mathbf{f}$ very even, there are two nilpotent orbits with Jordan type 
$\mathbf{f}$. Thus, in this case, $\mathrm{Exc}(\mu)$ 
contains at least $k$ irreducible divisors. 
On the other hand, for the Jacobson-Morozov parabolic subgroup $P$, 
$b_2(G/P) = k-1$ when $\g = sp(V)$, or $\g = so(V)$ with $\dim V$ odd. 
When $\g = so(V)$ and $\dim V$ is even, we must be careful; if the 
flag type of $P$ is of the form $(p_1, ..., p_{k-1}; 2; p_{k-1}, ..., p_1)$, 
$b_2(G/P) = k$. This happens when $\dim V = 4n + 2$ and 
$\mathbf{d} = [2^{2n}, 1^2]$ or when $\dim V = 8m+ 4n + 4$ and 
$\mathbf{d} = [4^{2m}, 3, 2^{2n}, 1]$. In the latter case, $\mathbf{d}$ has 
a minimal degeneration $\mathbf{d} \geq \mathbf{f}$ with 
$\mathbf{f} = [4^{2m}, 2^{2n+2}]$, which is very even.  
Note that $b_2(G/P)$ coincides with the relative Picard number $\rho$ 
of the Jacobson-Morozov resolution. 
By these observations, we know that $\mu$ has  
at least $\rho$ exceptional divisors   
except when $\g = so(4n+2)$, $n \geq 1$ and $\mathbf{d} = [2^{2n}, 1^2]$. 
Therefore, $\tilde{\mathcal{O}}_{\mathbf{d}}$ are {\bf Q}-factorial 
in these cases. By (1.3.2) they have terminal singularities. 
When $\g = so(4n+2)$, $n \geq 1$ and $\mathbf{d} = [2^{2n}, 1^2]$, 
$\mathcal{O}_{\mathbf{d}}$ is a Richardson orbit and the Springer 
map gives a small resolution of $\bar{\mathcal{O}}_{\mathbf{d}}$. 
Therefore, $\tilde{\mathcal{O}}_{\mathbf{d}}$ has non-{\bf Q}-factorial 
terminal singularities. 
\vspace{0.2cm}

(1.5) {\em Induced orbits in classical Lie algebras}: 
Let $\mathbf{d} = [d_1^{s_1}, ..., d_k^{s_k}]$ be an 
$\epsilon$-admissible partition of $m$. According as $\epsilon = 1$ or 
$\epsilon = -1$, we put $G = SO(m)$ or $G = Sp(m)$ respectively. 
Assume that $\mathbf{d}$ does 
not have full members. In other words, for some $p$, 
$d_p \geq d_{p+1} + 2$ or $d_k \geq 2$.  We put $r = \Sigma_{1 \le j \le p}{s_j}$. 
Then $\mathcal{O}_{\mathbf{d}}$ is an induced orbit (cf. [C-M], 7.3). More explicitly, 
there are a parabolic subgroup $Q$ of $G$ with (isotropic) flag type 
$(r, m-2r, r)$ with Levi decomposition 
$\q = \l \oplus \n$, and a nilpotent orbit $\mathcal{O}'$ of $\l$ 
such that $\mathcal{O}_{\mathbf{d}} = \mathrm{Ind}^{\g}_{\l}(\mathcal{O}')$. 
Here, $\l$ has a direct sum decomposition $\l = gl(r) \oplus \g'$, where 
$\g'$ is a simple Lie algebra of type $B_{(m-2r-1)/2}$ (resp. $D_{(m-2r)/2}$, 
resp. $C_{(m-2r)/2}$) when $\epsilon = 1$ and $m$ is odd (resp. 
$\epsilon = 1$ and $m$ is even, resp. $\epsilon = -1$). 
Moreover, $\mathcal{O}'$ is a nilpotent orbit of $\g'$ with Jordan type 
$[(d_1-2)^{s_1}, ..., (d_p-2)^{s_p}, d_{p+1}^{s_{p+1}}, ..., d_k^{s_k}]$.          
Let us consider the generalized Springer map 
$$ \nu: G \times^Q(n(\q) + \bar{\mathcal{O}}') \to \bar{\mathcal{O}}_{\mathbf{d}}$$ (cf. 
(1.2)). 
\vspace{0.2cm}

{\bf Lemma (1.5.1)}. {\em The map $\nu$ is birational. In other words, 
$\mathrm{deg}(\nu) = 1$.} 
\vspace{0.15cm}

{\em Proof}. We only discuss the case $G = Sp(m)$.   
It is enough to prove 
that $\nu^{-1}(z)$ is a point for $z \in \mathcal{O}_{\mathbf{d}}$. 
We prepare two kinds of skew-symmetric 
vector space $V_d$ ($d$: even), and $W_{2d}$ ($d$: odd) as follows. 
The $V_d$ is a $d$ dimensional vector space with 
a basis $e_1$, $e_2$, ..., $e_d$. 
The skew-symmetric form $\langle \;,\;\rangle$ is defined in such a way that 
$\langle e_i, e_j\rangle = 1$
(resp. $\langle e_i,e_j\rangle = -1$)  when $i + j = d+1$ and $i >j$ (resp. $i+j = d+1$ and $i < j$), and 
$\langle e_i, e_j\rangle = 0$ when $i + j \ne d+1$. 
Let $z_d$ be the endomorphism of $V_d$ defined by  the $d \times d$matrix 
$J_d$ with $J_d(i,i+1) = 1$ ($1 \leq i \leq d/2$), $J_d(i,i+1) = -1$ ($d/2 + 1 \leq i \leq d-1$) 
and otherwise $J_d(i,j) = 0$.   
The $W_{2d}$ is a $2d$ dimensional vector space with a basis 
$f_1$, ..., $f_d$, 
$f'_1$, ..., $f'_d$. 
The skew-symmetric form $\langle \;,\;\rangle$ is defined in such a way that $\langle f_i, f'_j\rangle = 1$
($\langle f'_i, f_j\rangle = -1$) when $i + j = d+1$, and otherwise $\langle f_i, f'_j\rangle = 0$, 
$\langle f'_i, f_j\rangle = 0$.  Let $z_{2d}$ be the endomorphism defined by the matrix 
$J'_{2d}$ with $J'_{2d}(i, i+1) = 1$ ($1 \le i \le d-1$), $J'_{2d}(d,d+1) = 0$, 
$J'_{2d}(i,i+1) = -1$ ($d+1 \le i \le 2d-1$) and otherwise $J'_{2d}(i,j) = 0$.   
Note that, in the partition $\mathbf{d}$,  $s_i$ is even if $d_i$ is odd.  
When $d_i$ is even, we put  $U_i := V_{d_i}^{\oplus s_i}$ and 
define $z_i \in \mathrm{End}(U_i)$ by $z_i = z_{d_i}^{\oplus s_i}$. 
When $d_i$ is odd, we put $U_i := W_{2d_i}^{\oplus s_i/2}$ and 
define $z_i \in \mathrm{End}(U_i)$ by $z_i = z_{2d_i}^{\oplus s_i/2}$. 
Let us consider the skew-symmetric vector space $V$ defined by   
$$ V := \oplus_{1 \le i \le k} U_i. $$
 Then we may assume that 
$z$ is an endomorphism of $V$ defined by $z = \oplus z_i$. 
Each $U_i$ has a filtration $0 \subset U_{i,1} \subset U_{i,2} 
\subset ... \subset U_{i, d_i} = U_i$ defined by $U_{i,j} := \mathrm{Im}(z_i^{d_i-j})$.  
By definition, $U_{i,d_i-1} = (U_{i,1})^{\perp}$.    
The problem is to find an $r$ dimensional isotropic subspace $F$ of $V$ in such 
a way that $z(F) = 0$ and 
 $z\vert_{(F^{\perp}/F)}$ has 
Jordan type $[(d_1-2)^{s_1}, ..., (d_p-2)^{s_p}, d_{p+1}^{s_{p+1}}, ..., d_k^{s_k}]$.    
We shall prove that $F = \oplus_{1 \le i \le p}U_{i,1}$. 
First note  that $F \subset \oplus_{1 \le i \le k}U_{i,1}$ since 
$z(F) = 0$. Next, one can check that 
$F^{\perp} \subset \oplus_{i \geq 2}U_i \oplus U_{1,d_1-1}$. 
In fact, if not, then one can find some $v \in F^{\perp}/F$ with 
$z^{d_1-2}(v) \ne 0$, which contradicts that $z\vert_{F^{\perp}/F}$ 
has Jordan type  $[(d_1-2)^{s_1}, ..., (d_p-2)^{s_p}, d_{p+1}^{s_{p+1}}, ..., d_k^{s_k}]$. 
By taking the dual with respect to the skew-symmetric form, one has 
$U_{1,1} \subset F$. We put $F_2: = F \cap \oplus_{2 \le i \le k}U_{i,1}$. 
Let $(F_2)^{\perp}$ be the orthogonal complement of $F_2$ in 
$\oplus_{2 \le i \le k}U_i$ (not in $V$).  Then $$F^{\perp}/F = 
(\bigoplus_{2 \le j \le d_1-1}U_{1,j}) \oplus (F_2)^{\perp}/F_2. $$ 
Then $x\vert_{(F_2)^{\perp}/F_2}$ has Jordan type 
$[(d_2-2)^{s_1}, ..., (d_p-2)^{s_p}, d_{p+1}^{s_{p+1}}, ..., d_k^{s_k}]$. 
We apply the same argument to $F_2$ to conclude that 
$U_{2,1} \subset F_2$. In particular, $U_{2,1} \subset F$. 
In this way,  we can prove inductively that $U_{i,1} \subset F$ 
for $i \leq p$. Since $\dim (\oplus_{1 \le i \le p}U_{i,1}) = r$, 
we have $F = \oplus_{1 \le i \le p}U_{i,1}$.  
\vspace{0.2cm}

{\bf Remark (1.5.2)}. A nilpotent orbit is called {\em rigid} if it is not induced from any 
other nilpotent orbit in a proper Levi subalgebra of $\g$. In a simple Lie algebra $\g$
of type B, C or D,  
$\mathcal{O}_{\mathbf d}$ is rigid if and only if $\mathbf{d}$ has full members and any 
odd (resp. even) member $d_i$ does not have multiplicity $s_i = 2$ when $\epsilon = 1$ 
(resp. $\epsilon = -1$) (cf. [CM]).  By Corollary (1.4.3), for a rigid orbit $\mathcal{O}_{\mathbf d}$, 
$\tilde{\mathcal O}_{\mathbf d}$ has only {\bf Q}-factorial terminal singularities. But, 
even if $\tilde{\mathcal O}_{\mathbf d}$ has only {\bf Q}-factorial terminal singularities, 
$\mathcal{O}_{\mathbf d}$ is not necessarily rigid.  
Assume that an odd (resp. even) member $d_p$ has multiplicity $2$ 
when $\epsilon = 1$ (resp. $\epsilon = -1$).  In this case, we have an induction of {\em another type}.  
Namely,  put $r = \Sigma_{1 \le j \le p-1}s_j + 1$. Then, there are a parabolic subgroup $Q$ of $G$ with 
(isotropic) flag type $(r, m-2r, r)$ with Levi decomposition 
$\q = \l \oplus \n$, and a nilpotent orbit $\mathcal{O}'$ in $\l$  
such that $\mathcal{O}_{\mathbf{d}} = \mathrm{Ind}^{\g}_{\l}(\mathcal{O}')$.
The orbit $\mathcal{O}'$ is contained in $\g'$, and its Jordan type is  
$[(d_1-2)^{s_1}, ..., (d_{p-1}-2)^{s_{p-1}}, (d_p-1)^2, d_{p+1}^{s_{p+1}}, ..., d_k^{s_k}]$.    
As explained above, if $\mathbf{d}$ has full members, but $\mathcal{O}_{\mathbf d}$ 
is not rigid, then $\mathcal{O}_{\mathbf d}$ has an induction of this type. 
But, for such an induction, the generalized Springer map $\nu$ is not birational.

\begin{center}
{\bf \S 2.} {\bf Main Results}
\end{center}

(2.1) Let $X$ be a complex algebraic variety with rational Gorenstein 
singularities. 
A partial resolution $f: Y \to X$ of $X$ is said to be a 
{\em $\mathbf{Q}$-factorial terminalization} of $X$ 
if $Y$ has only ${\mathbf Q}$-factorial 
terminal singularities and $f$ is a birational projective morphism 
such that $K_Y = f^*K_X$. In particular, when $Y$ is smooth, $f$ is called 
a crepant resolution. In general, $X$ has no crepant resolution; however, 
by [B-C-H-M], $X$ always has a $\mathbf{Q}$-factorial terminalization.  
But, in our case, the $\mathbf{Q}$-factorial terminalization can be 
constructed very explicitly without using the general theory in [B-C-H-M].   
\vspace{0.2cm}

{\bf Proposition (2.1.1)}. {\em Let $\mathcal{O}$ be a nilpotent orbit 
of a classical simple Lie algebra $\g$. 
Let $\tilde{\mathcal{O}}$ be the normalization of $\bar{\mathcal{O}}$. 
Then one of the following holds:}

(1) {\em $\tilde{\mathcal{O}}$ has $\mathbf{Q}$-factorial 
terminal singularities.} 
\vspace{0.15cm}

(2) {\em There are a parabolic subalgebra $\q$ of $\g$ with  
Levi decomposition $\q = \l \oplus \n$ and a nilpotent orbit $\mathcal{O}'$ 
of $\l$ such that (a): $\mathcal{O} = \mathrm{Ind}^{\g}_{\l}(\mathcal{O}')$ 
and (b): the normalization of $G \times^Q (n(\q) + \bar{\mathcal{O}'})$ is a  
$\mathbf{Q}$-factorial terminalization of $\tilde{\mathcal{O}}$ via 
the generalized Springer map.}   
  
{\em Proof}. When $\g$ is of type $A$, every $\tilde{\mathcal{O}}$ has a 
Springer resolution; hence (2) always holds.  
Let us consider the case $\g$ is of $B$, $C$ or $D$. 
Assume that (1) does not hold. Then, by (1.4.3),  the Jordan type $\mathbf{d}$ 
of $\mathcal{O}$ does not have full members except when 
$\g = so(4n+2)$, $n \geq 1$ and $\mathbf{d} = [2^{2n}, 1^2]$. 
In the exceptional case, $\mathcal{O}$ is a Richardson orbit and the Springer 
map gives a crepant resolution of $\tilde{\mathcal{O}}$; hence (2) holds. 
Now assume that $\mathbf{d}$ does not have full members. 
Then, by (1.5), $\mathcal{O}$ is an induced nilpotent orbit and there is a 
generalized Springer map 
$$\nu: G \times^Q(n(\q) + \bar{\mathcal{O}}') \to \bar{\mathcal{O}}.$$ 
This map is birational by (1.5.1). 
Let us consider the orbit $\mathcal{O}'$ instead of $\mathcal{O}$. 
If (1) holds for $\mathcal{O}'$, then $\nu$ induces a {\bf Q}-factorial 
terminalization of $\tilde{\mathcal O}$.  
If (1) does not hold for $\mathcal{O}'$, then 
$\mathcal{O}'$ is an induced orbit. 
By (1.2.3), one can replace $Q$ with a smaller parabolic subgroup $Q'$ 
in such a way that $\mathcal{O}$ is induced from $(Q', \mathcal{O}_2)$ for 
some nilpotent orbit $\mathcal{O}_2 \subset \l(Q')$. 
The generalized Springer map $\nu'$ for $(Q', \mathcal{O}_2)$ 
is factorized as 
$$ G \times^{Q'}(\n' + \bar{\mathcal{O}}_2) \to G \times^Q(\n + \bar{\mathcal{O}'}) \to \bar{\mathcal{O}}.$$ 
The second map is birational as explained above. 
The first map is locally obtained by a base change of 
the generalized Springer map 
$$L(Q) \times^{L(Q) \cap Q'}(\n(L(Q) \cap Q') + \bar{\mathcal{O}}_2) \to \bar{\mathcal{O}}'.$$ 
This map is birational by (1.5.1). Therefore, the first map is also birational, and 
$\nu'$ is birational. This induction step terminates and (2) finally holds. 
\vspace{0.2cm} 

(2.2) We shall next show that {\em every} $\mathbf{Q}$-factorial terminalization 
of $\tilde{\mathcal{O}}$ is of the form in Proposition (2.1.1) except 
when $\tilde{\mathcal{O}}$ itself has $\mathbf{Q}$-factorial terminal 
singularities. In order to do that, we need the following Proposition. 
\vspace{0,2cm}

{\bf Proposition (2.2.1)}. {\em Let $\mathcal{O}$ be a nilpotent orbit of 
a classical simple Lie algebra $\g$. Assume that a $\mathbf{Q}$-factorial terminalization 
of $\tilde{\mathcal O}$ is given by the normalization of 
$G \times^Q (n(\q) + \bar{\mathcal{O'}}))$ for some $(Q, \mathcal{O}')$ as in  
(2.1.1). Assume that $Q$ is a maximal parabolic subgroup of $G$ (i.e. 
$b_2(G/Q) = 1$), and this $\mathbf{Q}$-factorial terminalization is 
small.  
Then $Q$ is a parabolic subgroup corresponding to one of the following 
marked Dynkin diagrams and $\mathcal{O}' = 0$:}   
\vspace{0.15cm}

$A_{n-1}$ $(k < n/2)$

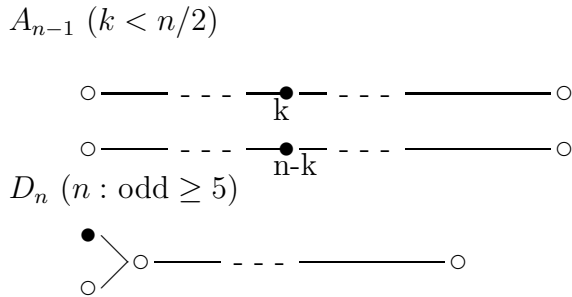
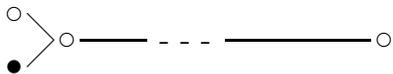
\begin{picture}(300,20)(0,0) 
\put(30,0){\circle{5}}\put(35,0){\line(1,0){25}} 
\put(65,-3.5){- - -}\put(90,0){\line(1,0){15}}
\put(105,0){\circle*{5}}\put(110,0){\line(1,0){10}}
\put(100,-10){k}\put(125,-3.5){- - -}\put(150,0)
{\line(1,0){55}}\put(210,0){\circle{5}}    
\end{picture}  

\begin{picture}(300,20)(0,0) 
\put(30,0){\circle{5}}\put(35,0){\line(1,0){25}} 
\put(65,-3.5){- - -}\put(90,0){\line(1,0){15}}
\put(105,0){\circle*{5}}\put(110,0){\line(1,0){10}}
\put(100,-10){n-k}\put(125,-3.5){- - -}\put(150,0)
{\line(1,0){55}}\put(210,0){\circle{5}}    
\end{picture}
\vspace{0.15cm}
 
$D_n$ $(n:$ $\mathrm{odd} \geq 5)$  

\begin{picture}(300,20)(0,0) 
\put(30,10){\circle*{5}}\put(30,-10){\circle{5}}
\put(35,10){\line(1,-1){10}}\put(35,-10){\line(1,1){10}}
\put(50,0){\circle{5}}\put(55,0){\line(1,0){25}}
\put(85,-3.5){- - -}\put(110,0){\line(1,0){55}}\put(170,0)
{\circle{5}}
\end{picture}
\vspace{0.7cm}

\begin{picture}(300,20)(0,0) 
\put(30,10){\circle{5}}\put(30,-10){\circle*{5}}
\put(35,10){\line(1,-1){10}}\put(35,-10){\line(1,1){10}}
\put(50,0){\circle{5}}\put(55,0){\line(1,0){25}}
\put(85,-3.5){- - -}\put(110,0){\line(1,0){55}}\put(170,0)
{\circle{5}}
\end{picture} 
\vspace{0.7cm}

{\em Proof}. Assume that $\g$ is of type $A$. Then every nilpotent 
orbit closure has a crepant resolution via the Springer map. 
Once $\tilde{\mathcal O}$ has a crepant resolution, 
every {\bf Q}-factorial terminalization is a crepant resolution (cf. [Na 2]).    
Then the claim follows from [Na], Proposition 5.1. Assume that $\g$ is of type $B$, $C$ 
or $D$. Since the {\bf Q}-factorial terminalization is small, 
$\tilde{\mathcal{O}}$ has terminal singularities. By (1.3.2) and 
(1.4.3), $\mathcal{O}$ is the nilpotent orbit in $so(4n+2)$, $n \geq 1$ 
with Jordan type $[2^{2n}, 1^2]$. In this case, $\mathcal{O}$ is a 
Richardson orbit and $\tilde{\mathcal{O}}$ has exactly two crepant resolutions 
via the Springer maps $G \times^Q n(\q)  \to \tilde{\mathcal{O}}$ 
corresponding to two marked Dynkin diagrams of type $D$ listed above (cf. [Na], p.92).    
Therefore, $\tilde{\mathcal{O}}$ has no other {\bf Q}-factorial terminalizations. 
\vspace{0.2cm}

The following is the main theorem: 
\vspace{0.2cm}

{\bf Theorem (2.2.2)}. {\em Let $\mathcal{O}$ be a nilpotent orbit 
of a classical simple Lie algebra $\g$. Then $\tilde{\mathcal{O}}$ 
always has a $\mathbf{Q}$-factorial terminalization. If  
$\tilde{\mathcal{O}}$ itself does not have $\mathbf{Q}$-factorial 
terminal singularities, then every $\mathbf{Q}$-factorial terminalization 
is given by the normalization of  
$G \times^Q (n(\q) + \bar{\mathcal{O'}}))$ in (2.1.1). 
Moreover, any two such $\mathbf{Q}$-factorial terminalizations 
are connected by a sequence of Mukai flops of type $A$ or $D$ defined 
in [Na], pp. 91, 92.}    

{\em Proof}. The first statement is nothing but (2.1.1). The proof of the 
second statement is quite similar to that of [Na], Theorem 6.1. 
Assume that $\tilde{\mathcal{O}}$ does not have $\mathbf{Q}$-factorial 
terminal singularities. Then, by (2.1.1), one can find a generalized 
Springer (birational) map 
$$\nu: G \times^Q(n(\q) + \bar{\mathcal{O}}') \to \bar{\mathcal{O}}.$$ 
Let $X_Q$ be the normalization of  $G \times^Q(n(\q) + \bar{\mathcal{O}}')$. 
Then $\nu$ induces a {\bf Q}-factorial terminalization 
$f: X_Q \to \tilde{\mathcal{O}}$. The relative nef cone $\overline{\mathrm{Amp}}(f)$ 
is a rational, simplicial, polyhedral cone of dimension $b_2(G/Q)$ (cf. (1.2.2) and 
[Na], Lemma 6.3).   
Each codimension one face $F$ of $\overline{\mathrm{Amp}}(f)$ corresponds to a 
birational contraction map $\phi_F : X_Q \to Y_Q$. 
The construction of $\phi_F$ is as follows. The parabolic subgroup $Q$ 
corresponds to a marked Dynkin diagram $D$. In this diagram, there are exactly 
$b_2(G/Q)$ marked vertexes. Choose a marked vertex $v$ from $D$. 
The choice of $v$ determines a codimension one face $F$ of $\overline{\mathrm{Amp}}(f)$.  
Let $D_v$ be the maximal, connected, single marked Dynkin subdiagram of $D$ 
which contains $v$.  Let $\bar{D}$ be the marked Dynkin diagram obtained from 
$D$ by erasing the marking of $v$. 
Let $\bar{Q}$ be the parabolic subgroup of $G$ corresponding to $\bar{D}$. 
Then, as in (1.2.2), we have a map  
$$\pi: G \times^Q (\n + \bar{\mathcal{O}}') \to G \times^{\bar{Q}} (\bar{\n} + \bar{\mathcal{O}_1}).$$  
Let $Y_Q$ be the normalization of $G \times^{\bar{Q}} (\bar{\n} + \bar{\mathcal{O}_1})$. 
Then $\pi$ induces a birational map $X_Q \to Y_Q$. This is the map $\phi_F$. 
Note that $\pi$ is locally obtained by a base change of the generalized Springer map 
$$L(\bar{Q}) \times^{L(\bar{Q}) \cap Q}
(\n(L(\bar{Q}) \cap Q) + \bar{\mathcal{O}'}) \to \bar{\mathcal{O}}_1.$$
Let $Z(\l(\q))$ (resp. $Z(\l(\bar{\q}))$) be the center of $\l(\q)$ (resp. $\l(\bar{\q})$). 
By the definition of $\bar{Q}$, the simple factors of $\l(\bar{\q})/Z(\l(\bar{\q}))$ are 
common to those of  $\l(\q)/Z(\l(\q))$ except one factor, say $\m$. 
Put  $\mathcal{O}'' := \mathcal{O}' \cap \m$.  
By (2.2.1), $\pi$ (or $\phi_F$) is a small birational map if and only if $\mathcal{O}'' = 0$ and $D_v$ 
is one of the single Dynkin diagrams listed in (2.2.1).  In this case, one can make a new 
marked Dynkin diagram$D'$ from $D$ by replacing $D_v$ by its dual $D^*_v$ (cf. [Na], Definition 1).   
Let $Q'$ be the parabolic subgroup of $G$ corresponding to $D'$. We may assume 
that $Q$ and $Q'$ are both contained in $\bar{Q}$. 
The Levi part of $Q'$ is conjugate to that of $Q$; hence there is a nilpotent orbit 
in $\l(\q')$ corresponding to $\mathcal{O}'$. We denote this orbit by the same 
$\mathcal{O}'$. Then $\mathcal{O}$ is induced from ($Q'$, $\mathcal{O}'$). 
As above, let $X_{Q'}$ be the normalization of 
$G \times^Q(n(\q') + \bar{\mathcal{O}}')$. Then we have a birational map 
$\phi'_F: X_{Q'} \to Y_Q$. The diagram $$X_Q \to Y_Q \leftarrow X_{Q'}$$ is a 
flop.  Assume that $g: X \to \tilde{\mathcal{O}}$ is a {\bf Q}-factorial terminalization. 
Then, the natural birational map $X --\to X_Q$ is an isomorphism in codimension 
one. Let $L$ be a $g$-ample line bundle on $X$ and let $L_0 \in \mathrm{Pic}(X_Q)$ 
be its proper transform of $L$ by this birational map. 
If $L_0$ is $f$-nef, then $X = X_Q$ and $f  = g$. Assume that $L_0$ is not 
$f$-nef. Then one can find a codimension one face $F$ of 
$\overline{\mathrm{Amp}}(f)$ which is negative with respect to $L_0$. 
Since $L_0$ is $f$-movable, the birational map 
$\phi_F: X_Q \to Y_Q$ is small. Then, as seen above, there is a new (small) birational map 
$\phi'_F: X_{Q'} \to Y_Q$.  Let $f': X_{Q'} \to \tilde{\mathcal O}$ be the composition 
of $\phi'_F$ with the map $Y_Q \to \tilde{\mathcal{O}}$. Then $f'$ is a {\bf Q}-factorial 
terminanization of $\tilde{\mathcal{O}}$. Replace $f$ by this $f'$ and repeat the same 
procedure; but this procedure ends in finite times (cf. [Na], Proof of Theorem 6.1 on 
pp. 104, 105).  More explicitly, there is a 
finite sequence of {\bf Q}-factorial terminalizations of $\tilde{\mathcal O}$:  
$$ X_0(:= X_Q) --\to  X_1(:= X_{Q'}) --\to .... --\to X_k (= X_{Q_k})$$ 
such that $L_k \in \mathrm{Pic}(X_k)$ is $f_k$-nef. This means that 
$X = X_{Q_k}$.
\vspace{0.2cm}

{\bf Example (2.3)}.  We put $G = SP(12)$. 
Let $\mathcal{O}$ be the nilpotent orbit in  
$sp(12)$ with Jordan type $[6, 3^2]$. Let $Q_1 \subset G$ be 
a parabolic subgroup with flag type $(3,6,3)$.  The Levi part $\l_1$ of 
$\q_1$ has a direct sum decomposition 
$$ {\mathfrak l}_1 = {\mathfrak g}l(3) \oplus C_3, $$ where $C_3$ is isomorphic to $sp(6)$ 
up to center.  Let $\mathcal{O}'$ be the nilpotent orbit in $sp(6)$ 
with Jordan type $[4, 1^2]$. 
Then $\mathcal{O} = \mathrm{Ind}^{sp(12)}_{\l_1}({\mathcal O}').$ 
Next consider the parabolic subgroup $Q_2 \subset SP(6)$ with flag 
type $(1,4,1)$.  The Levi part $\l_2$ of $\q_2$ has a direct sum 
decomposition $$ {\mathfrak l}_2 = {\mathfrak g}l(1) \oplus C_2 $$ where $C_2$ is isomorphic 
to $sp(4)$ up to center. Let $\mathcal{O}''$ be the nilpotent orbit 
in $sp(4)$ with Jordan type $[2, 1^2]$. Then   
$\mathcal{O}' = \mathrm{Ind}^{sp(6)}_{\l_2}({\mathcal O}'').$ 
One can take a parabolic subgroup $Q$ of $SP(12)$ with flag type $(3,1,4,1,3)$ 
in such a way that the Levi part $\l$ of $\q$ contains the nilpotent orbit 
$\mathcal{O}''$. Then $\mathcal{O}$ is the nilpotent orbit induced from 
$\mathcal{O}''$.  We shall illustrate the induction step above by 
$$ ([2, 1^2], sp(4)) \to ([4, 1^2], sp(6)) \to ([6, 3^2], sp(12)). $$  
Since $\tilde{O}''$ has only {\bf Q}-factorial terminal singularities, 
the {\bf Q}-factorial terminalization of $\tilde{\mathcal O}$ is given by the 
generalized Springer map 
$$ \nu: G \times^Q (n(\q) + \bar{\mathcal O}'') \to \bar{\mathcal O}. $$ 
The induction step is not unique; we have another induction step 
$$ ([2, 1^2], sp(4)) \to ([4,3^2], sp(10)) \to ([6,3^2], sp(12)). $$ 
By these inductions, we get another generalized Springer map 
$$ \nu': G \times^{Q'}(n(\q') + \bar{\mathcal O}'') \to \bar{\mathcal O}, $$ 
where $Q'$ is a parabolic subgroup of $G$ with flag type $(1,3,4,3,1)$. 
This gives another {\bf Q}-factorial terminalization of $\tilde{\mathcal O}$. 
The two {\bf Q}-factorial terminalizations of $\tilde{\mathcal O}$ are 
connected by a Mukai flop of type $A_3$.

\quad \\
\quad\\

Yoshinori Namikawa \\
Department of Mathematics, 
Graduate School of Science, Kyoto University, JAPAN \\
namikawa@math.kyoto-u.ac.jp


\begin{thebibliography}{} 

\bibitem[BCHM]{BCHM}  Birkar, C., Cascini, P., Hacon, C., McKernan, J.: 
Existence of minimal models for varieties of general type, math.AG/0610203

  
\bibitem[Be]{Be} Beauville,  A. : Symplectic singularities,  
Invent.  Math.  {\bf 139} (2000),  541-549 
 
\bibitem[C-G]{C-G} Chriss,  M. ,  Ginzburg,  V. : 
Representation theory and complex geometry,  
Progress in Math. ,  Birkhauser,  1997 

\bibitem[C-M]{C-M} Collingwood,  D. ,  McGovern,  W. : 
Nilpotent orbits in semi-simple Lie algebras,  
van Nostrand Reinhold,  Math.  Series,  1993 

\bibitem[Fu]{Fu} Fu,  B. : 
Symplectic resolutions 
for nilpotent orbits,  Invent.  Math.  {\bf 151}.  
(2003),  167-186

\bibitem[Hi]{Hi} Hinich,  V. : On the singularities of 
nilpotent orbits,  Israel J.  Math.  {\bf 73} (1991),  297-308 

\bibitem[Ka]{Ka} Kawamata, Y.: Crepant blowing-up of 3-dimensional 
canonical singularities and its application to degenerations of 
algebraic surfaces, Ann. Math. {\bf 127} (1988), 93-163 

\bibitem[K-P]{K-P} Kraft, H., Procesi, C.: 
On the geometry of conjugacy classes in classical groups, 
Comment. Math. Helv. {\bf 57}. (1982), 539-602

\bibitem[L-S]{L-S} Lusztig, G., Spaltenstein, N.: 
Induced unipotent classes, J. London Math. Soc. 
{\bf 19}. (1979), 41-52

\bibitem[Na]{Na} Namikawa,  Y. : Birational geometry 
of symplectic resolutions of nilpotent orbits,  
Advances Studies in Pure Mathematics {\bf 45},  
(2006),  Moduli Spaces and Arithmetic Geometry (Kyoto,  2004),  
pp.  75-116 

\bibitem[Na 2]{Na 2} Namikawa, Y.: Flops and Poisson deformations 
of symplectic varieties, Publ. RIMS {\bf 44}. No 2.  (2008), 259-314 

\bibitem[Na 3]{Na 3} Namikawa, Y.: Extensions of 2-forms and 
symplectic varieties, J. Reine Angew. Math. {\bf 539} (2001), 123-147

\bibitem[Na 4]{Na 4} Namikawa, Y.: Birational geometry and deformations of 
nilpotent orbits, Duke Math. J. {\bf 143} (2008), 375-405

\bibitem[Pa]{Pa} Panyushev,  D. I. : Rationality of singularities and the 
Gorenstein property of nilpotent orbits,  Funct.  Anal.  Appl.  {\bf 25} 
(1991),  225-226

\bibitem[Ri]{Ri} Richardson,  R. : Conjugacy classes in parabolic subgroups 
of semi-simple algebraic groups,  Bull.  London Math.  Soc.  {\bf 6} 
(1974),  21-24



 
\end{thebibliography}
\end{document}